\documentclass[11pt,a4paper]{article}

\usepackage{latexsym,amssymb, amsthm}
\usepackage[centertags]{amsmath}
\usepackage{epsfig,pifont}
\usepackage[inactive]{srcltx}  

\usepackage{pgf,pgfarrows,pgfnodes,pgfautomata,pgfheaps,pgfshade}
\usepackage{graphics,helvet,times,mathptm,epic,eepic}

\usepackage{color}

\usepackage{newlfont}

\usepackage{amsfonts,bbm}

\newenvironment{keywords}{ \noindent {\small\bf Key Words}:}{ }

\def\bd{\begin{description}}
\def\ed{\end{description}}

\def\beq{\begin{equation}}
\def\eeq{\end{equation}}
\def\bea{\begin{eqnarray}}
\def\eea{\end{eqnarray}}
\def\beas{\begin{eqnarray*}}
\def\eeas{\end{eqnarray*}}

\theoremstyle{remark}

\begin{document}

\title{\textbf{\textsc{Evaluating the exact infinitesimal values of area of   Sierpinski's carpet and volume of   Menger's sponge}}}

\newcommand{\nms}{\small}

\author{ {   \bf Yaroslav D. Sergeyev\footnote{Distinguished Full Professor, he works also  at the N.I.~Lobatchevsky State University,
  Nizhni Novgorod, Russia (Professor, part-time contract) and at    the Institute of High Performance
  Computing and Networking of the National Research Council of Italy (affiliated researcher).    } } \\ \\ [-2pt]
      \nms Dipartimento di Elettronica, Informatica e Sistemistica,\\[-4pt]
       \nms   Universit\`a della Calabria,  Via P. Bucci, Cubo 42-C,\\[-4pt]
       \nms 87030 Rende (CS)  -- Italy\\
       \nms http://wwwinfo.deis.unical.it/$\sim$yaro\\[-4pt]
         \nms {\tt  yaro@si.deis.unical.it }
}

\date{}

\maketitle

\vspace{-1cm}

\begin{abstract}

Very often traditional approaches     studying dynamics of
self-similarity processes are not able to give their quantitative
characteristics at infinity and, as a consequence, use limits to
overcome this difficulty. For example, it is well know that the
limit area of Sierpinski's carpet and volume of   Menger's sponge
are equal to zero. It is shown  in this paper  that recently
introduced infinite and infinitesimal numbers allow us to use
exact expressions instead of limits and to calculate exact
infinitesimal values of areas and volumes at various points at
infinity  even if the chosen moment of the observation  is
infinitely faraway on the time axis from the starting point. It is
interesting that traditional results that can be obtained without
the usage of infinite and infinitesimal numbers can be produced
just as finite approximations of the new ones. The importance of
the possibility to have this kind of quantitative characteristics
for \textit{\textbf{E}}-Infinity theory is emphasized.
 \end{abstract}


\section{Introduction}
\label{se1}

Appearance of  new powerful approaches modelling  the spacetime by
fractals (see
\cite{El_Naschie_2004,El_Naschie_2005,El_Naschie,He_2005,He_2006,Iovane_Part_1,Iovane_Part_2,Sidharth}
and references given therein) urges  development of  adequate
mathematical tools allowing one to study fractal objects
quantitatively after $n$ steps executed in a fractal process for
both finite and  \textit{infinite} $n$. Very well developed
theories of fractals (see e.g.,
\cite{El_Naschie,Falconer,Hastings_Sugihara,fractals,fractals2}
and references given therein) allow us to give certain numerical
answers to questions regarding fractals (calculation of, e.g.,
their length, area or volume) only for   finite values of $n$. The
same questions very often remain without any answer when we
consider an infinite number of steps because  the traditional
mathematics
  (both standard and non-standard versions of Analysis) can
speak only about limit fractal objects and the required values
tend to zero or infinity.

\begin{figure}[t]
  \begin{center}
    \epsfig{ figure = 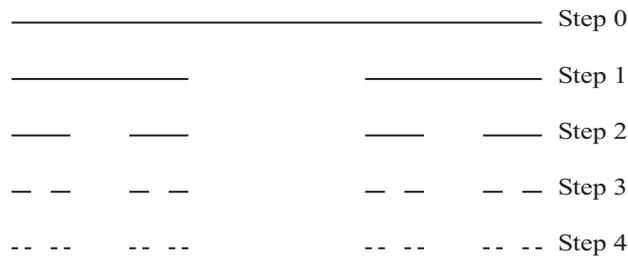, width = 3.2in, height = 1.3in,  silent = yes }
    \caption{Cantor's set.}
 \label{Fractals_6}
  \end{center}
\end{figure}

 Let us consider, for example, the famous Cantor's set
(see Fig.~\ref{Fractals_6}).   If a finite number of steps, $n$,
has been done in   construction of Cantor's set, then we are able
to describe numerically the set being the result of this
operation. It will have $2^n$ intervals having the length
$\frac{1}{3^n}$ each. Obviously, the set obtained after $n+1$
iterations will be different and we also are able to measure the
lengths of the intervals forming the second set.   It will have
$2^{n+1}$ intervals having the length $\frac{1}{3^{n+1}}$ each.
The situation changes drastically in the limit because we are not
able to distinguish results of $n$ and $n+1$ steps of the
construction if $n$ is infinite.

We also are not able to distinguish at infinity the results of the
following two processes that both use Cantor's construction but
start from different positions. The first one is the usual
Cantor's set and it starts from the interval $[0,1]$, the second
starts from the couple of intervals $[0,\frac{1}{3}]$ and
$[\frac{2}{3},1]$. In spite of the fact that for any given finite
number of steps, $n$, the results of the constructions will be
different for these two processes we have no traditional
mathematical tools allowing us to distinguish them at infinity.

Recently a new applied point of view on infinite and infinitesimal
numbers (that is not related to the non-standard analysis of
Robinson \cite{Robinson}) has been introduced in
\cite{Sergeyev,Poland,chaos,informatica,Korea,Dif_Calculus}. The
new approach evolves Cantor's ideas about existence of different
infinite numbers. It describes infinite and infinitesimal numbers
that are in accordance with the principle `The part is less than
the whole'  and gives a possibility to work with finite, infinite,
and infinitesimal quantities \textit{numerically} by using a new
kind of a computer -- Infinity Computer -- introduced in
\cite{Sergeyev_patent,www,Poland}. A detailed analysis of Cantor's
set by using these new tools has been done in \cite{chaos}. A
comprehensive introduction to the new approach and examples of its
usage can be found in \cite{www,informatica}).

In this paper, we show how the  new computational tools can be
used to calculate the exact infinitesimal values of area of
Sierpinski's carpet and volume of   Menger's sponge. Calculation
of numerical characteristics of fractal objects of this kind is
particularly important in the context of
\textit{\textbf{E}}-Infinity theory (see
\cite{El_Naschie_2004,El_Naschie_2005,El_Naschie,El_Naschie_2008,
He_2005,He_2006,Iovane_Part_1,Iovane_Part_2,Iovane_Part_3,Sidharth}).

\section{Physical methodology in Mathematics}
\label{se2}

In Physics,   researchers use tools to describe the object of
their study and the used instrument influences results of
observations and restricts possibilities of observation of the
object. Thus, there exists the philosophical triad -- researcher,
object of investigation, and tools used to observe the object. A
new applied approach to infinity proposed in
\cite{Sergeyev,informatica,Korea,Dif_Calculus}  emphasizes
existence of this triad in Mathematics, as well. Mathematical
languages (in particular, numeral systems\footnote{ We remind that
\textit{numeral}  is a symbol or group of symbols that represents
a \textit{number}. For example, the symbols `10', `ten', and `X'
are different numerals, but they all represent the same number.})
are among the tools used by mathematicians to observe and to
describe mathematical objects. Very often difficulties that we
find solving mathematical problems are related not to their nature
but to inadequate mathematical languages used to solve them.  The
new approach is based on the following methodological postulates
(see \cite{Sergeyev,informatica,Korea}).

\textbf{Postulate 1.} \textit{There exist infinite and
infinitesimal objects but   human beings and machines are able to
execute only a finite number of operations.}

\textbf{Postulate 2.} \textit{We shall not   tell \textbf{what
are} the mathematical objects we deal with; we just shall
construct more powerful tools that will allow us to improve our
capacities to observe and to describe properties of mathematical
objects.}

\textbf{Postulate 3.} \textit{The principle `The part is less than
the whole' is applied to all numbers (finite, infinite, and
infinitesimal) and to all sets and processes (finite and
infinite).}

Due to this declared applied statement, such concepts as
bijection, numerable and continuum sets, cardinal and ordinal
numbers cannot be used in this paper because they belong to the
theories working with different assumptions. As a consequence, the
new approach   is different also with respect to the non-standard
analysis introduced in \cite{Robinson} and built using Cantor's
ideas. It is important to emphasize that our point of view on
axiomatic systems is also more applied than the traditional one.
Due to Postulate~2,    mathematical objects are not define by
axiomatic systems that just determine formal rules for operating
with certain numerals reflecting some properties of the studied
mathematical objects.

Due to Postulate 3, infinite and infinitesimal numbers are managed
in the same manner as we are used to deal with finite ones. This
Postulate in our opinion very well reflects organization of the
world around us but   in many traditional infinity theories   it
is true only for finite numbers. Due to Postulate~3, the
traditional point of view on infinity accepting such results as
$\infty + 1= \infty$ are substituted in  such a way   that $\infty
+ 1 > \infty$.

This methodological program has been realized in
\cite{Sergeyev,informatica,Korea} where a new powerful numeral
system has been developed. This system gives   a possibility to
execute \textit{numerical} computations not only with finite
numbers but also with infinite and infinitesimal ones in
accordance with Postulates 1--3. The main idea consists of
measuring  infinite and infinitesimal quantities   by different
(infinite, finite, and infinitesimal) units of measure.

A new infinite unit of measure   has been introduced for this
purpose in \cite{Sergeyev,informatica,Korea} in accordance with
Postulates 1--3 as the number of elements of the set $\mathbb{N}$
of natural numbers. It is expressed by a new numeral \ding{172}
called \textit{grossone}. It is necessary to emphasize immediately
that the infinite number \ding{172} is not either Cantor's
$\aleph_0$ or $\omega$. Particularly, it has both cardinal and
ordinal properties as usual finite natural numbers. Formally,
grossone is introduced as a new number by describing its
properties postulated by the \textit{Infinite Unit Axiom} (IUA)
(see \cite{Sergeyev,informatica,Korea}). This axiom is added to
axioms for real numbers similarly to addition of the axiom
determining zero to axioms of natural numbers when integer numbers
are introduced.

Inasmuch as it has been postulated that grossone is a number,  all
other axioms for numbers hold for it, too. Particularly,
associative and commutative properties of multiplication and
addition, distributive property of multiplication over addition,
existence of   inverse  elements with respect to addition and
multiplication hold for grossone as for finite numbers.
Introduction of grossone gives a possibility to work with finite,
infinite, and infinitesimal quantities \textit{numerically} by
using a new kind of a computer -- the Infinity Computer --
introduced in \cite{Sergeyev_patent,www,Poland}.

\section{Evaluating infinitesimal areas and volumes}
\label{se3}

Grossone can be successfully used for various purposes related to
infinite and infinitesimal objects, particularly, for measuring
infinite sets (see \cite{Sergeyev,www,chaos,informatica}).  For
instance, the number of elements of a set
$B=\mathbb{N}\backslash\{b\}$, $b \in \mathbb{N}$,
 is equal to $\mbox{\ding{172}}-1$ and the
number of elements of a set $A=\mathbb{N}\cup\{a\}$, where $a
\notin \mathbb{N}$, is equal to $\mbox{\ding{172}}+1$. Note that
due to Postulate~3 and definition of \ding{172}, the number
$\mbox{\ding{172}}+1$ is not natural (see \cite{informatica} for a
detailed discussion). Analogously, $\mbox{\ding{172}}^3$ is the
number of elements of the set
  $W$, where
 \beq
   W  =
\{ (a_1, a_2, a_3)  : a_1 \in   \mathbb{N}, a_2 \in \mathbb{N},
a_3 \in   \mathbb{N} \}.
 \label{Ext_1}
\eeq Grossone can   be used to indicate positions of elements in
infinite sequences, too. Let us investigate this issue in detail.
First of all, we remind that the traditional definition of the
infinite sequence is: \textit{`An infinite sequence $\{a_n\}, a_n
\in A $ for all $n \in \mathbb{N},$ is a function having as the
domain the set of natural numbers, $\mathbb{N}$, and as the
codomain  a set $A$'}. We have postulated    that the set
$\mathbb{N}$ has \ding{172} elements. Thus, due to the sequence
definition given above, any sequence having $\mathbb{N}$ as the
domain  has \ding{172} elements. Traditionally, the notion of a
subsequence is introduced as a sequence from which some of its
elements have been cancelled. This definition gives sequences
having the number of members (finite or infinite) less than
grossone. Consider, for instance, these two infinite sequences
 \beq
\underbrace{1,2,3,4,\hspace{1mm}  \ldots \hspace{1mm}
\mbox{\ding{172}}-2,\hspace{1mm}
 \mbox{\ding{172}}-1,
\mbox{\ding{172}}}_{\mbox{\ding{172} elements}}, \hspace{5mm}
\underbrace{3,4,5,6, \hspace{1mm}  \ldots \hspace{1mm}
\mbox{\ding{172}}-2,\hspace{1mm}
 \mbox{\ding{172}}-1}_{\mbox{\ding{172}-3 elements}}.  \label{Ext_2}
 \eeq
Thus, the number of elements of any   sequence (finite or
infinite) is less or equal to~\ding{172}. One of the immediate
consequences of the understanding of this result is that any
sequential process can have at maximum \ding{172} steps. It is
very important to notice a deep relation of this observation to
the Axiom of Choice. Since \textit{any} sequential process can
have at maximum \ding{172} elements, this is true for the process
of choice, as well. Therefore, it is not possible to choose in a
sequence more than \ding{172} elements from a set. For example, if
we consider the following set
\[
D = \{1,2,3,4,\hspace{1mm}  \ldots \hspace{1mm}
\mbox{\ding{172}}-2,\hspace{1mm}
 \mbox{\ding{172}}-1,
\mbox{\ding{172}}, \mbox{\ding{172}}+1, \mbox{\ding{172}}+2,
\ldots \mbox{\ding{172}}^3-2, \mbox{\ding{172}}^3-1,
\mbox{\ding{172}}^3 \}
\]
and   a process of the sequential choice of elements  from this
set in the order of increase of their values, it is possible to
arrive at maximum to~\ding{172} starting from 1
 \beq
\underbrace{1,2,3,4,\hspace{1mm}  \ldots \hspace{1mm}
\mbox{\ding{172}}-2,\hspace{1mm}
 \mbox{\ding{172}}-1,
\mbox{\ding{172}}}_{\mbox{\ding{172} steps} },
\mbox{\ding{172}}+1, \mbox{\ding{172}}+2,  \ldots
\mbox{\ding{172}}^3-2, \mbox{\ding{172}}^3-1, \mbox{\ding{172}}^3
\label{Ext_3}
 \eeq
executing so  \ding{172} steps. Starting from 3 it is possible
to arrive at maximum to $\mbox{\ding{172}}+2$
 \beq
1,2,\underbrace{3,4,\hspace{1mm}  \ldots \hspace{1mm}
\mbox{\ding{172}}-2,\hspace{1mm}
 \mbox{\ding{172}}-1,
\mbox{\ding{172}},  \mbox{\ding{172}}+1,
\mbox{\ding{172}}+2}_{\mbox{\ding{172} steps}},
\mbox{\ding{172}}+3,
  \ldots \mbox{\ding{172}}^3-2,
\mbox{\ding{172}}^3-1, \mbox{\ding{172}}^3. \label{Ext_4}
 \eeq
Of course, due to Postulate 1, we are able to observe only a
finite number among \ding{172} members of these processes.  In
addition, it depends on the chosen numeral system, $\mathcal{S}$,
which numbers
  we can observe because $\mathcal{S}$ should have numerals able to represent
the required numbers.

\begin{figure}
  \begin{center}
\epsfig{ figure = 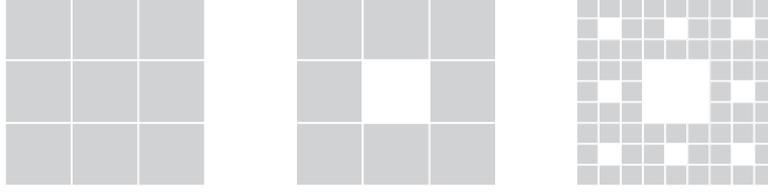, width = 4in, height = 1in,
silent = yes }
    \caption{Sierpinski's carpet.} \label{Extensions_1}
  \end{center}
  \end{figure}

Another important observation consists of the fact that numeral
systems including \ding{172} allow  us to observe the starting and
the \textit{ending} elements of infinite processes, if the
respective elements are expressible in these numeral systems. This
fact is very important in connection with fractals because it
allows us to distinguish different fractal objects after an
infinite number of steps of their construction. This observation
concludes preliminary results and allows us to start to study
Sierpinski's carpet (see Fig.~\ref{Extensions_1}).

Since   construction of Sierpinski's carpet is a process,  it
cannot contain more then \ding{172} steps (see discussion related
to   (\ref{Ext_2})-(\ref{Ext_4})). Thus, if at iteration $n=1$ the
process starts from the first from the left grey square having the
length of a side equal to 1, at iteration $n \ge 1$ Sierpinski's
carpet consists   of $N_n=8^{n-1}$ grey boxes. The length of a
side of a hole, $L_n$, is equal to $3^{-(n-1)}$  and  the
 area of all grey boxes     is equal to
\beq
 A_n = L_n^{2}\cdot N_n = \left(\frac{8}{9}\right)^{n-1}.
 \label{Ext_5}
\eeq
 Therefore, if   \ding{172} steps have been executed, the area
$A_{\tiny{\mbox{\ding{172}}}} =
\left(\frac{8}{9}\right)^{\tiny{\mbox{\ding{172}}}-1}$ and if
\ding{172}-9 steps have been executed, the area
$A_{\tiny{\mbox{\ding{172}}}-9} =
\left(\frac{8}{9}\right)^{\tiny{\mbox{\ding{172}}}-10}$. It is
worthwhile  to notice   that (again due to the limitation
illustrated by the example (\ref{Ext_2})-(\ref{Ext_4})) it is not
possible to count one by one all the   boxes at Sierpinski's
carpet if their number is superior to \ding{172}. For instance, at
iteration \ding{172}   it has $8^{\tiny{\mbox{\ding{172}}}-1}$
boxes and they cannot be counted sequentially because
$8^{\tiny{\mbox{\ding{172}}}-1}
> \mbox{\ding{172}}$ and any process (including that of the
sequential counting) cannot have more that \ding{172} steps.

Thus, we are able now to distinguish different Sierpinski's
carpets at different points at infinity and to calculate the
respective areas that are expressed in infinitesimals. Moreover,
we can do it also when we change the starting element of the
construction. For instance, if at iteration $n=1$   the process
starts from the second from the left object consisting of 8 grey
squares having the length of a side equal to $\frac{1}{3}$, at
iteration $n \ge 1$ Sierpinski's carpet consists of
$N_{2,n}=8^{n}$ grey boxes, where the subscript indicates the
numbers of the starting and the ending points of the process. The
length of a side of a hole $L_{2,n}=3^{-n}$  and  the
 area of all grey boxes     is equal to
\beq
 A_{2,n} = L_{2,n}^{2}\cdot N_{2,n} = \left(\frac{8}{9}\right)^{n}.
 \label{Ext_6}
\eeq
 Then, if   \ding{172} steps have been executed, the area $A_{2,\tiny{\mbox{\ding{172}}}}
= \left(\frac{8}{9}\right)^{\tiny{\mbox{\ding{172}}}}$ and if
  \ding{172}-9 steps have been executed, the area $A_{2,\tiny{\mbox{\ding{172}}}-9} =
\left(\frac{8}{9}\right)^{\tiny{\mbox{\ding{172}}}-9}$. Obviously,
formulae (\ref{Ext_5}), (\ref{Ext_6}) can be easily generalized to
the   case where $N_{k,n}=8^{n+k-2},$ $ L_{k,n}=3^{-(n+k-2)}$ and
 \beq
 A_{k,n} = L_{k,n}^{2}\cdot N_{k,n} =
 \left(\frac{8}{9}\right)^{n+k-2}, \hspace{5mm} 1 \le k \le n \le \mbox{\ding{172}} +k -1,
 \label{Ext_7}
\eeq
 and  (\ref{Ext_7}) can be used
   both for finite and infinite values of $k$ and $n$.

\begin{figure}
  \begin{center}
        \epsfig{ figure = 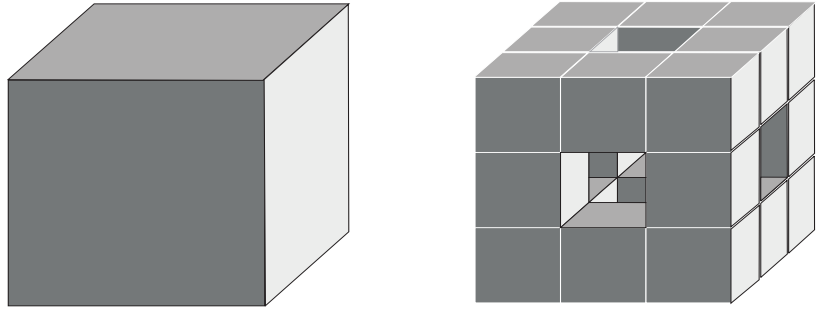, width = 4.2in, height = 1.6in,  silent = yes }
    \caption{Menger's sponge.}
 \label{Extensions_2}
  \end{center}
  \end{figure}

We conclude this paper by evaluating the volume of Menger's sponge
by a complete analogy to (\ref{Ext_5})--(\ref{Ext_7}). Let $n$ be
the number of iterations. Then, starting from the left cube at
iteration $n=1$ the number of grey cubes is $N_{1,n}=20^{n-1}, n
\ge 1$. The length of a side of a hole $L_{1,n}$ is equal to
$3^{-(n-1)}$ and the
 volume of all grey cubes at the $n$-th iteration is equal to
\beq
 V_{1,n} = L_{1,n}^{3}\cdot N_{1,n} = \left(\frac{20}{27}\right)^{n-1}.
 \label{Ext_8}
\eeq
 For $n =\mbox{\ding{172}}$ we have $V_{1,\tiny{\mbox{\ding{172}}}} =
\left(\frac{20}{27}\right)^{\tiny{\mbox{\ding{172}}}-1}$ and for
$n =\mbox{\ding{172}}-1$ it follows
$V_{1,\tiny{\mbox{\ding{172}}}-1} =
\left(\frac{20}{27}\right)^{\tiny{\mbox{\ding{172}}}-2}$. Finally,
the general formula of the volume of Menger's sponge is
 \beq
 V_{k,n} = L_{k,n}^{3}\cdot N_{k,n} = \left(\frac{20}{27}\right)^{n+k-2}, \hspace{5mm} 1 \le k \le n \le \mbox{\ding{172}} +k -1.
 \label{Ext_9}
\eeq

\section{Conclusion}

Very often traditional approaches     studying dynamics of
self-similarity processes are not able to give their quantitative
characteristics at infinity and, as a consequence, use limits to
overcome this difficulty. For example, it is well know that the
limit area of Sierpinski's carpet and volume of   Menger's sponge
are equal to zero.

In this paper, it has been shown that  infinite and infinitesimal
numbers  introduced in
\cite{Sergeyev,Poland,chaos,informatica,Korea,Dif_Calculus} allow
us to obtain exact numerical results instead of limits and to
calculate exact infinitesimal values of areas and volumes of
various fractal objects  at different points at infinity. In fact,
the possibility to express explicitly various infinite numbers
allows us to indicate the final elements not only for finite but
for infinite processes, as well. As a result, we can calculate the
required areas and volumes at the chosen   moment even if this
moment is infinitely faraway on the time axis from the starting
point. It is interesting that traditional results that can be
obtained without the usage of infinite and infinitesimal numbers
can be produced just as finite approximations of the new ones.

It has been shown the importance of the starting conditions in
fractal processes. Suppose that   there are two identical fractal
processes starting from different starting structures in such a
way that the starting structure of the second process is a result
of $k$ steps of the first one. Then, after the same number of
steps, $n$, the two process lead to different results for both
finite and infinite values of $n$. Thus, this fact being obvious
for finite values of $n$ holds for infinite values of $n$, too.
Finally, the importance of the possibility to have this kind of
quantitative characteristics for \textit{\textbf{E}}-Infinity
theory (see
\cite{El_Naschie_2004,El_Naschie_2005,El_Naschie,El_Naschie_2008,
He_2005,He_2006,Iovane_Part_1,Iovane_Part_2,Iovane_Part_3,Sidharth})
 has been emphasized in the paper.

\bibliographystyle{plain}
\bibliography{XBib_Extensions}
\end{document}